\documentclass[12pt,a4paper,notitlepage]{article}
\usepackage{amsfonts}
\begin{document}

\title{\textbf{ON QUANTUM GROUPS CO-REPRESENTATIONS}}
\author{\textbf{H. Montani}$^{*}$\textbf{\ and R. Trinchero}$^{+}$ \\
%EndAName
\medskip\ \\
\textit{Centro At\'{o}mico Bariloche,}\\
\textit{\ Comisi\'{o}n Nacional \ de Energ\'{\i}a At\'{o}mica,}\\
\textit{\ and}\\
\textit{\ Instituto Balseiro, }\\
\textit{\ Universidad Nacional de Cuyo.}\\
\textit{\ 8400 - S. C. de Bariloche, Rio Negro.}\\
\textit{\ Argentina.}}
\date{}
\maketitle

\begin{abstract}
We carry out a generalization of quantum group co-representations in order
to encode in this structure those cases where non-commutativity between
endomorphism matrix entries and quantum space coordinates happens.
\end{abstract}

\vspace*{\fill}

*e-mail: montani@cab.cnea.edu.ar

+e-mail: trincher@cab.cnea.edu.ar

\newpage\

\section{Introduction}

Quantum Groups arise as the abstract structure underlying the symmetries of
integrable systems in (1+1) dimensions\cite{fad}. There, the theory of
quantum inverse scattering give rise to some deformed algebraic structures
which were first explained by Drinfeld as deformations of classical Lie
algebras\cite{drinf}\cite{jim}. Also, an analog structure was obtained by
Woronowicz in the context of non-commutative $C^{*}$-algebras\cite{wor}.
There is a third approach, due to Manin, where Quantum Groups are
interpreted as the endomorphisms of certain non commutative algebraic
varieties defined by quadratic algebras, called quantum linear spaces (QLS)%
\cite{man}\cite{kass}. In this work, we are mainly involved with this last
approach.

From the algebraic point of view, quantum groups are Hopf algebras and the
relation with the endomorphism algebra of QLS come from their
co-representations on tensor product spaces. It is worth remarking that the
usual construction of the co-action on the tensor product space implies the
commutativity between the matrix elements of a representation of the
endomorphism and the coordinates of the QLS.

In the present work, we are interested in relaxing this last statement thus
admitting non-commutative relations between endomorphism matrix entries and
quantum linear space coordinates. We carry out this task modifying the
definition of co-action on tensor product spaces. The question is: at what
extent can this be done without spoling out the bialgebra and its
co-representation structure? We answer this question, finding the conditions
under which the general framework of quantum groups and quantum linear space
still holds. This construction introduces a new family of deformation
parameters in the endomorphism bialgebra of the QLS, as it is shown in the
quantum plane example where a two parameter $M_{q,p}(2)$ is obtained. Also,
we find a non central object playing the rol of $(q,p)$-deformed determinant
which allows to obtain a $GL_{q,p}(2)$.

We present a brief description of the co-representations of bialgebras in
the second section, developing our approach to modified co-representations
in the third one and, finally, we present the quantum plane example in the
last section.

\section{Quantum algebras and co-representations}

Let $V$ be a vector space of dimension $n$, $\{e_i\}$ a basis for $V$ and $%
H_o$ the trivial bialgebra of functions over $GL(n,\mathbf{C})$. This
bialgebra is freely generated by the identity and the coordinates functions $%
T_i^j$, in the basis $\{e_i\}$, defined by

\[
T{_i^j:GL(n,\mathbf{C})\longrightarrow }\mathbf{C}
\]

\[
T{_i^j:g\longmapsto g_i^j}
\]
for $g\in GL(n,\mathbf{C})$.

The co-product and co-unit are given by

\begin{equation}
\Delta T_i^j=T_i^k\otimes T_k^j  \label{comult}
\end{equation}

\begin{equation}
\varepsilon (T_i^j)=\delta _j^i  \label{coun}
\end{equation}

From now on, sumation over repeated index is assumed.

The comodule $(\delta ,V)$, with

\[
{\delta :V\longrightarrow H_o{\otimes }V}
\]

\begin{equation}
\delta ({e_i)=T_i^j{\otimes }e_j}  \label{coac}
\end{equation}
provides a representation of $GL(n,\mathbf{C})$ in $V$, through the $g_i^j$
in the basis $b$ of $V$. It has the co-associativity property and preserves
the co-unit, which is expressed by the commutativity of the diagrams

\begin{center}
\begin{equation}
\begin{tabular}{ccc}
$V$ & $\stackrel{\delta _V}{\longrightarrow }$ & $H_o\otimes V$ \\
$^{\delta _V}\downarrow $ &  & $\downarrow ^{I_{H_o}\otimes \delta _V}$ \\
$H_o\otimes V$ & $\stackrel{\Delta \otimes I_V}{\longrightarrow }$ & $%
H_o\otimes H_o\otimes V$%
\end{tabular}
\label{coass}
\end{equation}

\begin{equation}
\begin{tabular}{ccccc}
&  & $H_o\otimes V$ &  &  \\
$^{\delta _V}$ & $\nearrow $ &  & $\searrow $ & $^{\epsilon \otimes I_V}$ \\
$V$ &  & $\stackrel{I_V}{\longrightarrow }$ &  & $V$%
\end{tabular}
\label{counit}
\end{equation}
\end{center}

Building up co-representations for objects with more structure than $%
V^{\otimes }$, as quadratic algebras for example, requires some extra
conditions that we sketch below.

Let $A$ denote the quadratic algebra generated by the ideal $I(\frak{B}),$
where $\frak{B}:V\otimes V\longrightarrow V\otimes V$, then

\begin{equation}
A(\frak{B})=\frac{V^{\otimes }}{I(\frak{B})}  \label{alq}
\end{equation}
$V^{\otimes }$ is tensor algebra on $V.$ In general we consider $\frak{B}$
with the form
\begin{equation}
\frak{B}=(I_{V\otimes V}-B)  \label{bb}
\end{equation}

\begin{equation}
e_ie_j-B_{ij}^{kl}e_ke_l  \label{eebb}
\end{equation}

In order to built up a co-module structure on $A(\frak{B})$, one needs to
work on the tensor algebra $V^{\otimes }.$ This is obtained with the natural
extension of the co-action on $V\otimes V$ :

\begin{equation}
\delta _{V\otimes V}=(m\otimes I_V)(I_{H_o}\otimes \tau \otimes I_V)(\delta
_V\otimes \delta _V)  \label{coacvv}
\end{equation}
where $m:H_o\otimes H_o\longrightarrow H_o$ is the product in the bialgebra $%
H_o$ .It is worth remarking that the flip operator $\tau :V\otimes
H_o\longrightarrow H_o\otimes V$ will render a commutative product between
the ${T_i^j}$ and $e_k,$ as it is assumed in the construction of the quantum
plane and $GL_q(2)$\cite{man}

The main condition for preserving the quadratic algebra structure is the
''commutativity'' between $\frak{B}$ and $\delta _{V\otimes V}:$

\begin{center}
\begin{equation}
\begin{tabular}{ccc}
$V\otimes V$ & $\stackrel{\frak{B}}{\longrightarrow }$ & $V\otimes V$ \\
$^{\delta _{V\otimes V}:}\downarrow $ &  & $\downarrow ^{\delta _{V\otimes
V}:}$ \\
$H\otimes V\otimes V$ & $\stackrel{I_H\otimes \frak{B}}{\longrightarrow }$ &
$H\otimes V\otimes V$%
\end{tabular}
\label{bddb}
\end{equation}
\end{center}

This graph is satisfied if $H$ is the bialgebra arising from the quotient of
the free algebra generated by the objects ${T_i^j}$ and the ideal $I(\frak{B}%
,H_o)$ generated by

\begin{equation}
B_{ij}^{kl}T{_k^rT_l^s-T_i^kT_j^l}B_{kl}^{rs}  \label{btt}
\end{equation}
i.e.,

\[
H=\frac{H_o}{I(\frak{B},H_o)}
\]

Since $I(\frak{B},H_o)$ is a co-ideal with relation to $\Delta $ then $H$
becomes a bialgebra. The equation (\ref{btt}) is a central object in the so
called FRT construction\cite{frt}. In this way, $A(\frak{B})$ becomes in a $%
H $-algebra comodule.

\section{Generalized co-representations}

The main aim of this section is to build up the mathematical framework
encoding the situation in which entries of the endomorphism matrix may not
commute with the coordinates of the quantum linear space defined in (\ref
{alq}). We will reach it by means of a modification in the co-representation
theory. As described in the previous section, supplying the quantum linear
space with a $H$-comodule structure requires a good definition of a
co-action on $V\otimes V$.

As we saw, the standard definition of $\delta _{V\otimes V}$, eq.(\ref
{coacvv}), provides both $V\otimes V$ with a $H_o$-comodule structure and $A=%
\frac{V^{\otimes }}{I(\frak{B})}$ with an $H$-comodule structure.

We will show that one can replace the flip map $\tau $ by a more general one
$\gamma $ in such a way that there exist a $H_\gamma $-comodule structure,
for some $H_\gamma $ to be constructed, provided $\gamma $ satisfy some
requirements. Then, let $\gamma $ be defined by

\begin{equation}
\gamma :V\otimes H\longrightarrow H\otimes V  \label{g}
\end{equation}

\begin{equation}
\gamma (e_iT{_j^k)=\gamma }_{ijn}^{klm}{T}_l^n{e}_m  \label{egt}
\end{equation}

So, we propose the co-action on tensor product space as being

\begin{equation}
\delta _{V\otimes V}^\gamma =(m\otimes I_V)(I_{H_\gamma }\otimes \gamma
\otimes I_V)(\delta _V\otimes \delta _V)  \label{coacg}
\end{equation}

In order to have an $H_\gamma $-comodule structure, the co-associativity and
co-unity properties, diagrams (\ref{coass}) and (\ref{counit}), must be
satisfied together with requirement that $\frak{B}$ be a co-module map, (\ref
{bddb}).

Recalling the bijection between comodules and multiplicative matrices\cite
{man}, let us consider the multiplicative matrix $M$ in $V\otimes V$ with
coefficients in $H_o$ corresponding to the comodule $\delta _{V\otimes
V}^\gamma $ ,i.e.,
\[
\delta _{V\otimes V}^\gamma \equiv M\in End(V\otimes V,H_o)
\]

\begin{equation}
\delta _{V\otimes V}^\gamma (e_i\otimes e_j)=M_{ij}^{rs}\otimes e_r\otimes
e_s  \label{coacg1}
\end{equation}
hence $M$ is

\begin{equation}
M_{ij}^{kl}={\gamma }_{mjn}^{lpk}{T}_i^m{T}_p^n.  \label{M}
\end{equation}
Also, we adopt the following convention: let $A_{ij}^{kl}$ and $D_{ij}^{kl}$
any pair of four-tensors, then we write

\begin{equation}
(A\times B)_{ij}^{rs}=A_{ij}^{kl}\times D_{kl}^{rs}  \label{prod}
\end{equation}
where $\times $ stands for any kind of product (tensor, algebraic,etc.), and
sum over repeated index is assumed.

We start analyzing the diagram (\ref{bddb}), which means that the co-module $%
\delta _{V\otimes V}^\gamma $ map preserves the quantum linear space $A(%
\frak{B})$

\begin{equation}
A(\frak{B})=\frac{V^{\otimes }}{I(\frak{B})},  \label{alq}
\end{equation}
that it is expressed by the commutation relation

\[
(I_{H_\gamma }\otimes \frak{B})\circ \delta _{V\otimes V}^\gamma =\delta
_{V\otimes V}^\gamma \circ \frak{B}
\]

From this equation we get the condition

\[
\frak{B}M=M\frak{B}
\]
meaning that $H_\gamma $ must be defined as being

\[
H_\gamma =\frac{H_o}{I(\frak{B}M-M\frak{B})}
\]

Here $I(\frak{B}M-M\frak{B})$ is the ideal generated by the quadratic
relation $\frak{B}M-M\frak{B}.$

A necessary and sufficient condition for $H_\gamma $ to be a bialgebra is
that $I(\frak{B}M-M\frak{B})$ be a co-ideal,i.e.,

\begin{equation}
\Delta I\subset I\otimes H+H\otimes I  \label{coid}
\end{equation}

Here, it must be taken into account that in order to have a $H_\gamma $%
-comodule structure on $A(\frak{B})$, $\delta _{V\otimes V}^\gamma $ must
fulfill the co-associativity and co-unity properties, diagrams (\ref{coass})
and (\ref{counit}),

\[
(\Delta \otimes I_{V\otimes V})\circ \delta _{V\otimes V}^\gamma
=(I_{H_\gamma }\otimes \delta _{V\otimes V}^\gamma )\circ \delta _{V\otimes
V}^\gamma
\]

\[
(\epsilon \otimes I_{V\otimes V})\circ \delta _{V\otimes V}^\gamma
=I_{V\otimes V},
\]

Again, because of the bijection between all the structures of left co-module
on $V=\mathbf{C}^n$ and the multiplicative matrix $\mathbf{M(}n,H\gamma )$%
\cite{man}, if $M\in H_\gamma $ satisfy
\begin{equation}
\Delta M=M\otimes M  \label{dmmm}
\end{equation}

\begin{equation}
\epsilon (M)=I  \label{em}
\end{equation}
automatically both properties are satisfied.

Then, coming back to the co-ideal question, eq. (\ref{coid}), and assuming
that this last condition holds, one gets

\begin{eqnarray*}
\Delta {(}\frak{B}M-M\frak{B}) &=&\frak{B}\Delta M-M\Delta \frak{B} \\
&=&\frak{B}(M\otimes M)-(M\otimes M)\frak{B} \\
&=&{(}\frak{B}M-M\frak{B})\otimes M-M\otimes {(}\frak{B}M-M\frak{B})
\end{eqnarray*}
and

\[
\epsilon {(}\frak{B}M-M\frak{B})=\frak{B}\epsilon (M)-\epsilon (M)\frak{B}=0
\]
hence $H_\gamma $ is a bialgebra.

This results may be resumed into the following assertion: if there exist $%
\gamma $ satisfying $\frak{B}M-M\frak{B}$ , $\Delta M=M\otimes M$ and $%
\varepsilon (M)=I_{V\otimes V},$ for $M\in H_\gamma =\frac{H_o}{I(\frak{B}M-M%
\frak{B})}$, then $H_\gamma $ is a bialgebra and $\delta _{V\otimes
V}^\gamma $ renders $A(\frak{B})$ into a $H_\gamma $-comodule.

Now, we consider the antipode. If one can define an antipode , i.e., a map $%
S:H_\gamma \longrightarrow H_\gamma $ such that the following property holds

\begin{equation}
m\circ (S\otimes I_{H_\gamma })\circ \Delta =m\circ (I_{H_\gamma }\otimes
S)\circ \Delta =\eta \circ \epsilon  \label{antipoda}
\end{equation}
at the level of the multiplicative matrices representing the bialgebra $%
H_\gamma $ in $A(\frak{B})$ via a comodule $(A(\frak{B}),\delta ),$ then we
can assert that there exists a representation of the resulting Hopf algebra,
given by the comodule $(A(\frak{B})\otimes A(\frak{B}),\delta _{V\otimes
V}^\gamma )$. In fact recalling that,

\begin{equation}
S(M_{ij}^{kl})={\gamma }_{mjn}^{lpk}S({T}_p^n)S({T}_i^m).  \label{antt}
\end{equation}
it is simple to see that,
\[
m\circ (S\otimes I_{H_\gamma })\circ \Delta M=m\circ (I_{H_\gamma }\otimes
S)\circ \Delta M=\eta \circ \epsilon (M)
\]
and using (\ref{dmmm}) ensures that $S(M)$ is the inverse of $M,$

\[
MS(M)=S(M)M=I
\]

So, we may affirm that the construction presented above still holds if $%
H_\gamma $ is a Hopf algebra, then encoding the situation of non-commutative
co-representations in the framework of Quantum Groups.

In the next section, we describe an explicit example enjoying all these
properties presented above, namely a $p$-deformed version of the
endomorphism of the quantum plane.

\section{The Quantum Plane}

Let the quantum plane $A_q^{2\mid 0}$ described by

\begin{equation}
e_1e_2=qe_2e_1  \label{qpq}
\end{equation}

In the basis $\{e_1\otimes e_{1,}e_1\otimes e_2,e_2\otimes e_1,e_2\otimes \}$%
, this relation can be expressed by means the quadratic form $B$ as

\begin{equation}
B=\left[
\begin{array}{cccc}
1 & 0 & 0 & 0 \\
0 & 0 & q & 0 \\
0 & q & 1-q^2 & 0 \\
0 & 0 & 0 & 1
\end{array}
\right]  \label{bqp}
\end{equation}
which is a solution of the Yang-Baxter equation

\begin{equation}
B_{12}B_{23}B_{12}=B_{23}B_{12}B_{23}  \label{yb}
\end{equation}

It is worth remarking that the following construction leads to the same
structure for other choice of $B$, as the symmetric and idempotent $%
B^{\prime },$

\begin{equation}
B^{\prime }=\frac 1{q+q^{-1}}\left[
\begin{array}{cccc}
q+q^{-1} & 0 & 0 & 0 \\
0 & q-q^{-1} & 2 & 0 \\
0 & 2 & q^{-1}-q & 0 \\
0 & 0 & 0 & q+q^{-1}
\end{array}
\right]  \label{bpq1}
\end{equation}

This $B^{\prime }$ is not a solution of Yang-Baxter equation but, in the
Manin construction for pseudo-symmetric quantum space\cite{man}, it enables
to characterize all the endomorphism of the quantum plane by the relation $%
B^{\prime }M-MB^{\prime }$ as the only solution to the master relation $%
(I-B^{\prime })M(I+B^{\prime }).$

The endomorphism matrix $T$ is

\begin{equation}
T=\left[
\begin{array}{cc}
a & b \\
c & d
\end{array}
\right]  \label{tqp}
\end{equation}

Then, proposing the simplest form for $\gamma ,$

\begin{equation}
{\gamma }_{ijn}^{klm}=\delta _i^m\delta _j^l\delta _n^kg(i,j,k)  \label{ddg}
\end{equation}
we get $M,$

\begin{equation}
M_{ij}^{kl}=g(k,j,l)T_i^kT_j^l  \label{mqp}
\end{equation}

Thus, solving the pair of conditions

\[
\Delta M=M\otimes M
\]

\[
BM=MB
\]
we obtain the following $g(i,j,k)$

\begin{equation}
\begin{tabular}{l}
$g(i,j,j)=1$ \\
\\
$g(i,1,2)=p$ \\
\\
$g(i,2,1)=p^{-1}$%
\end{tabular}
\label{gqp}
\end{equation}
where $p\in \mathbf{C}$ is new parameter. Introducing this $g$ in $\gamma ,$
and substituting in the relation $BM=MB$ one gets

\begin{equation}
\begin{array}{r}
ac-pqca=0 \\
\\
ab-p^{-1}qba=0 \\
\\
bc-p^2cb=0 \\
\\
cd-p^{-1}qdc=0 \\
\\
bd-pqdb=0 \\
\\
ad-da+p(q^{-1}-q)cb=0
\end{array}
\label{bmqp}
\end{equation}
which define a two parametric $M_{q,p}(2)$ as the quotient algebra

\begin{equation}
M_{q,p}(2)=\frac{k[T_i^j]}{I(BM-MB)}  \label{mq2}
\end{equation}
for $i$ and $j$ from 1to 2.

Also, this $g$ give rise to the following relations between matrix entries
and the coordinates of the quantum plane,

\begin{equation}
\begin{array}{r}
e_ia-ae_i=0 \\
\\
e_id-de_i=0 \\
\\
e_ib-pbe_i=0 \\
\\
e_ic-p^{-1}ce_i=0
\end{array}
\label{emqp}
\end{equation}

The same construction can be carried out on the Grassmannian plane $%
A_q^{0\mid 2}$ given by the quadratic relations

\begin{equation}
e_1e_2=-q^{-1}e_2e_1  \label{gqpq}
\end{equation}

\begin{equation}
e_i^2=0  \label{nh}
\end{equation}
leading to same set of relations (\ref{bmqp}) and (\ref{emqp}), so $%
M_{q,p}(2)$ is also the bialgebra of endomorphism of $A_q^{0\mid 2}$.

It is also possible to define an antipode on $M_{q,p}(2)$. First at all, we
need to define something like a determinant, but there is no such an object
in the center of $M_{q,p}(2).$ However we may built up a $\det_{q,p}$ such
that it has homogenous commutation relation with the generator of $%
M_{q,p}(2),$ i.e., $\{a,b,c,d\}.$ Then, we define $D=\det\nolimits_{q,p}\in $
$\frac{k[a,b,c,d]}{I(BM-MB)},$

\begin{equation}
D=\det\nolimits_{q,p}=ad-p^{-1}qbc  \label{detq}
\end{equation}
that satisfy the following commutation relations

\begin{center}
\begin{equation}
\begin{array}{r}
Da-aD=0 \\
\\
Db-p^{-2}bD=0 \\
\\
Dc-p^2cD=0 \\
\\
Dd-dD=0
\end{array}
\label{dabcd}
\end{equation}
\end{center}

With these properties, and assuming that $D$ is an invertible element of $%
M_{q,p}(2)$ we define the antipode $S$

\begin{equation}
S\left[
\begin{array}{cc}
a & b \\
c & d
\end{array}
\right] =D^{-1}\left[
\begin{array}{cc}
d & -(pq)^{-1}b \\
-pqc & a
\end{array}
\right]   \label{antip}
\end{equation}
and, consequently, we have the Hopf algebra $GL_{q,p}(2).$

\section{Concluding Remarks}

We have introduced a new ingredient in the theory of representations of
Quantum (semi)Groups admitting non-commutativity between endomorphism matrix
entries and quantum space coordinates. This feature give rise to an extra
deformation of all the involved structures. In fact, we have shown that it
is possible to introduce the map $\gamma :V\otimes H_\gamma \longrightarrow
H_\gamma \otimes V$ without spoiling out the Hopf algebra and co-module
structures provided that $\gamma $ turns $H_\gamma $ into a Quantum Matrix
Group. This was explicitly shown in the Quantum Plane example, where we
obtained a two parameter deformation $M_{q,p}(2)$ and also we were able to
construct a $H_\gamma $ valued function $\det_{q,p}$, out of the centre of
the algebra, thus making possible to get $GL_{q,p}(2).$

\section{Acknowledgment}

We are greatly indebted to M. L. Bruschi for enlightening discussions. Also,
the authors thank CONICET-Argentina for finantial support.

\newpage\

\end{document}